\documentclass[11pt]{article}
\usepackage{amssymb}
\usepackage{color}
\usepackage{graphicx}
\usepackage{float}

\RequirePackage{amsthm, amsmath, amsfonts, amssymb,epsf}%, graphicx}
\usepackage{dsfont}
\newcommand{\be}{\begin{equation}}
\newcommand{\ee}{\end{equation}}

\newcommand{\ep}{\varepsilon}
\newcommand{\var}{  \hbox{\rm var} }
\newcommand{\sP}{{\mathcal P}}
\newcommand{\sD}{{\mathcal D}}
\newcommand{\R}{{\mathbb R}}
\newcommand{\Z}{{\mathbb Z}}
\newcommand{\half}{ {\scriptstyle{\frac{1}{2}} } }
\newcommand{\quarter}{ {\scriptstyle{\frac{1}{4}} } }

\newtheorem{proposition}{Proposition}
\newtheorem{theorem}{Theorem}
\newtheorem{defin}{Definition}

\begin{document}
\date{\today}
\title{The least favorable noise}
\author{Philip A. Ernst\footnote{Department of Statistics, Rice University}, Abram M. Kagan\footnote{Department of Mathematics, University of Maryland}, and L.C.G. Rogers\footnote{Statistical Laboratory, University of Cambridge}}
\maketitle
\begin{center}
\small \textit{We dedicate this work to our colleague, mentor, and friend, Professor Larry Shepp (1936–2013)}.
\end{center}

\begin{abstract}
\noindent
Suppose that a random variable $X$ of interest is observed perturbed by independent additive noise $Y$. This paper concerns the ``the least favorable perturbation'' $\hat Y_\ep$, which maximizes the prediction error $E(X-E(X|X+Y))^2$ in the class of $Y$ with $ \var (Y)\leq \ep$. We find a characterization of the answer to this question, and show by example that it can be surprisingly complicated. However, in the special case where $X$ is infinitely divisible, the solution is complete and simple.
We also explore the conjecture that noisier $Y$ makes prediction worse.
\end{abstract}
{\bf Keywords}: Least favorable perturbation. self-decomposable random variable, infinitely divisible distributions.\\\\
{\bf MSC 2010 Subject Codes}: Primary: 60E07, 60E10 Secondary: 60E05
%%%%%%%%%%%%%%%%%%%%%%%%%%%%%%%%%%%%%%%%%%%%%%%%%%%%%%%%%%%%%%%%%%%%%%

\section{Introduction.}
Suppose that on a probability space we observe $X+Y$, where $X$ and $Y$ are independent random variables, $X$ being a square-integrable random variable of interest, and $Y$ being an additive noise perturbation.
The prediction error
\begin{equation}
E\{X-E(X|X+Y)\}^2 ={\rm var}X- {\rm var}E(X|X+Y)
\label{prederr}
\end{equation}
depends on $Y$ of course, and thus a natural question is `What would be the {\em worst} noise $Y$ we could add to $X$?'  In other words, given the law of $X$, how would we choose the law of $Y$ to maximize the prediction error in equation \eqref{prederr}, or equivalently, how would we find
\begin{equation}
\inf_Y \var E(X|X+Y) \;?
\end{equation}
Since the mean of $E[X|X+Y]$ is fixed and equal to $EX$, an equivalent question is to choose the law of $Y$ so as to achieve
\begin{equation}
\inf_Y E \bigl\lbrace \; E[X|X+Y]^2 \; \bigr\rbrace.
\label{min2}
\end{equation}
If we think  of what happens when $Y = \lambda Z$, where $Z \sim N(0,1)$, we quickly realize that as $\lambda \rightarrow \infty$  we have
\begin{equation}
E[X|X+\lambda Z] \rightarrow EX \qquad \hbox{\rm a.s.,}
\label{normnoise}
\end{equation}
so that the minimization in \eqref{min2} has a trivial solution unless we bound the variance of $Y$. So we will focus on the problem
\begin{equation}
\inf\bigl\lbrace \; E \bigl( \; E[X|X+Y]^2 \; \bigr):
\; \var(Y) \leq \varepsilon \; \bigr\rbrace,
\label{min3}
\end{equation}
where $\varepsilon >0$ is given.  We then have a number of questions:
\begin{description}
\item[Question 1:]Can we find an explicit solution to \eqref{min3}?
\item[Question 2:] Can we characterize the solution to \eqref{min3}?
\item[Question 3:] Are there situations with explicit solutions?
\item[Question 4:]Does more noise mean worse prediction?
\end{description}
\indent \indent The fact that we asked Question 2 means that the answer to Question 1 has to be `No'; however, the answer to Question 2 is `Yes', and we deal with this in Section \ref{S1}. The answer to Question 3 is also `Yes', as we show in Section \ref{S2}; if the law of $X$ is infinitely divisible, then we can find  the minimizing $Y$. Simple examples show that the answer to Question 4 is `No', but if $Y$ is self-decomposable we have a partial result in this direction; see Section \ref{S2a}.  In Section \ref{S3}, we present an analysis of the case where $X$ is binomial and $Y$ is integer-valued, and we give a number of numerical examples which point to the diversity and complexity of the solutions in general. \\
\indent  We conclude with some brief remarks about the broader literature. The spirit of this work is most closely aligned with the lines of inquiry in \cite{Bryc0,Bryc,Dembo}. We also note that the focus of the present work largely moves in the opposite direction of stochastic filtering, in which one (usually) seeks to get as close as possible to $X$ (see, e.g., \cite{Crisan}, and references therein). This being said, the answer to Question 4 should be of interest to those in stochastic filtering.

\section{Characterizing the solution.}\label{S1}
Firstly, we observe that the objective to be minimized, 
\begin{equation}
 \var \bigl( \; E[X|X+Y] \; \bigr)
 \label{obj1}
\end{equation}
is unaltered if we shift $X$ or $Y$ by a constant, so we may and shall assume that the means of $X$ and $Y$ are set to be zero, unless otherwise stated.\\
\indent If $f$ is the density of $X$ and $g$ is the density of $Y$, then the objective \eqref{min2} is to minimize $\Phi(g)$, where
\begin{equation*}
\Phi(g) \equiv \int \; \lambda_*(s)^2\, (f*g)(s) \; ds,
\label{obj2}
\end{equation*}
and where
\begin{equation*}
\lambda_*(s) \equiv \frac{    \int xf(x)g(s-x) \; dx }{(f*g)(s)}\;\;
\label{ladef}
\end{equation*}
is $E[X|X+Y=s]$. We notice firstly that
\begin{equation*}
\Phi(g) =\sup_\lambda \int \bigl\{\; 2 \lambda(s)\, \int xf(x)g(s-x) \; dx
-\lambda(s)^2 (f*g)(s) \; \bigr\} \; ds,
\label{phi2}
\end{equation*}
which tells us in particular that $\Phi$ is a convex function. We aim to minimize $\Phi(g)$ over feasible $g$, that is, $g$ in the set $\sP$ of primal-feasible functions:
\begin{equation}
\sP = \{ g\geq 0, \int g(y)dy = 1, \int yg(y) dy = 0, \int y^2 g(y) dy \leq \varepsilon\}.
\end{equation}
Writing $\sigma^2$ for the second moment of $X$, we then have that for any $\alpha, \beta \in \R$ and $\gamma \geq 0$
\begin{eqnarray}
\inf_{g \in \sP} \Phi(g) &\geq & \inf_{g \in \sP, z \geq 0}
\bigl\{  \Phi(g) + \alpha (\int g dy -1) + \beta\int y g dy
+\gamma(\int y^2 g dy + z - \varepsilon) \bigr\}
\nonumber\\
&\geq&\inf_{g \geq 0, z \geq 0}
\bigl\{  \Phi(g) + \alpha (\int g dy -1) + \beta\int y g dy
+\gamma(\int y^2 g dy + z - \varepsilon) \bigr\}
\nonumber \\
&\geq&\inf_{g \geq 0, z \geq 0} \sup_\lambda \biggl[
\int \bigl\{\; 2 \lambda(s)\, \int xf(x)g(s-x) \; dx
-\lambda(s)^2 (f*g)(s) \; \bigr\} \; ds
\nonumber \\
&& \qquad + \alpha (\int g dy -1) + \beta\int y g dy
+\gamma(\int y^2 g dy + z - \varepsilon)
\biggr]
\nonumber\\
&\geq& \sup_\lambda \inf_{g \geq 0, z \geq 0}
\biggl[
\int \bigl\{\; 2 \lambda(s)\, \int xf(x)g(s-x) \; dx
-\lambda(s)^2 (f*g)(s) \; \bigr\} \; ds
\nonumber \\
&& \qquad + \alpha (\int g dy -1) + \beta\int y g dy
+\gamma(\int y^2 g dy + z - \varepsilon)
\biggr]
\nonumber\\  %%%%%%%%%%%%%%%%%%%%%%%%%%%%%%%%%%%%%%%%%%%%%%
&=& \sup_\lambda \inf_{g \geq 0, z \geq 0}
\biggl[
\int g(y)\bigl\{\alpha + \beta y + \gamma y^2
- \int f(s-y)(s-y-\lambda(s))^2
 \; ds 
)    
\bigr\}\; dy 
\nonumber\\
&&\qquad\qquad+\sigma^2-\alpha + \gamma (z-\varepsilon)
\biggr].
\label{end}
\end{eqnarray}
From this, we deduce that 
\begin{equation}
\inf_{g \in \sP} \Phi(g)\geq
\sup_{(\alpha, \beta,\gamma, \lambda) \in \sD}\bigl[ \; 
\sigma^2 -\alpha - \gamma\varepsilon
\; \bigr],
\label{primal_dual}
\end{equation}
where $\sD$ is the space of dual-feasible variables $(\alpha,
 \beta,\gamma, \lambda) $ satisfying $\gamma \geq 0$ and the condition
 \begin{equation}
 0 \leq \alpha + \beta y + \gamma y^2
- \int f(s-y)(s-y-\lambda(s))^2
 \; ds  \qquad\forall y.
 \label{df1}
 \end{equation}
The inequality in \eqref{primal_dual} is a primal-dual inequality familiar
from constrained optimization problems. We expect that under technical conditions
it is possible to prove that the inequality is in fact an equality, but 
we avoid attempting to prove this. We do so because establishing this (if true) 
does not help us to {\em identify} an optimal solution in any particular 
example; to do that we will have to exploit the special features of the solution, 
and by so doing we will be able to pass directly to a proof of optimality. We now proceed to do so.

\begin{theorem}\label{thm1}
Suppose that $g_* \in \sP$ and $(\alpha_*, \beta_*, \gamma_*, \lambda_*) \in \sD$
satisfy the complementary slackness conditions
\begin{eqnarray}
0 &\equiv& \biggl[\alpha_* + \beta_* y + \gamma_* y^2
- \int f(s-y)(s-y-\lambda_*(s))^2
 \; ds \biggr] \; g_*(y) 
 \label{cs1}
 \\
 0 &=& \gamma_* z_*,
 \label{cs2}
\end{eqnarray}
where $z_* = \varepsilon - \int y^2 g_*(y) \; dy $;  and that
\begin{equation}
\lambda_*(s) = \frac{    \int xf(x)g_*(s-x) \; dx }{(f*g_*)(s)}.
\label{ladef2}
\end{equation}
Then $g_*$ is optimal.
\end{theorem}

\noindent{\sc Proof.} Consider $q \equiv \sigma^2 - \alpha_* - \gamma_* \varepsilon$, which is a lower bound for the right-hand side of \eqref{primal_dual}, since $(\alpha_*, \beta_*, \gamma_*, \lambda_*) \in \sD$. Now we return to \eqref{end} and work back through the steps, putting $g_*$ for $g$ and $\lambda_*$ for $\lambda$, ignoring the $\sup$ and $\inf$ everywhere. Because of the conditions in \eqref{cs1} and \eqref{cs2}, the value we start from at \eqref{primal_dual} is $q$. At every step, we have equality, so we end up with $\Phi(g_*) = q$. Since $g_* \in \sP$, $g_*$ is optimal. This concludes the proof.

\hfill$\square$

It might appear that the conditions of Theorem \ref{thm1} are too complicated to verify in practice, but upon inspection of \eqref{cs1} we realize that {\em if for some $a, b \in \R$}
\begin{equation}
\lambda_*(s) = a + bs,
\label{good}
\end{equation}
{\em there may be a chance}.  Indeed, if we continue to assume that $EX=EY=0$, then the condition 
\begin{equation}
\alpha_* + \beta_* y + \gamma_* y^2 = \int f(x) (x-\lambda_*(x+y))^2 dx \qquad\forall y,
\label{cs1b}
\end{equation}
combined with \eqref{good} implies that $a=0$, $\beta_*=0$, $b^2 = \gamma_*$, and $\alpha_* = (1-b)^2 \sigma^2$.

\section{Explicitly soluble situations.}\label{S2}
If we took $X$, $Y$ to be independent with the {\em same} distribution, then it is obvious that 
\begin{equation}
E[X|X+Y] = (X+Y)/2.
\label{IID}
\end{equation}
\indent Let us now apply Theorem \ref{thm1} to this situation, taking $\lambda_*(s) = s/2$, $\alpha_* = \sigma^2/4$, $\beta_*=0$, and $\gamma_*=1/4$, and $g = f$. If the bound on the variance of $Y$ is $\varepsilon = \sigma^2 \equiv \var(X)$, then $g=f$ is primal-feasible, $(\alpha_*, \beta_*, \gamma_*, \lambda_*)$ is dual-feasible, and the complementary slackness conditions \eqref{cs1b} and \eqref{cs2} hold. Hence by Theorem \ref{thm1} the law which minimizes $\var E[X|X+Y]$ subject to the bound $\var(Y) \leq \var(X)$ is $g=f$. The lower bound from \eqref{primal_dual} is seen to be $\sigma^2/2$, which is indeed the variance of $(X+Y)/2$.\\
\indent By similar reasoning, it is straightforward to see that if $X = \xi_1 + \ldots + \xi_n$, where the $\xi_j$ are IID with zero mean and common variance $\sigma^2$, and where we bound $\var(Y) \leq m \sigma^2$, then the optimal law of $Y$ is given by $Y = \xi_1+ \ldots + \xi_m$. But this result now points towards a wider result for infinitely divisible distributions, which we state as Proposition \ref{Prop1} below.

\begin{proposition} \label{Prop1}
Suppose that $(Z_t)_{t \geq 0}$ is a zero-mean square-integrable L\'evy process, with $EZ_t = t$. Suppose that  $X \sim Z_t$ for some fixed $t>0$. Then the minimum in \eqref{min3} is achieved when $Y \sim Z_\varepsilon$.
\end{proposition}

\bigbreak\noindent
{\sc Proof.} If we let $Y = Z_\varepsilon$, then $E[X|X+Y] = t(X+Y)/(t+\varepsilon)$, so by setting $\lambda_*(s) = bs$ with $b=t/(t+\varepsilon)$ we ensure that \eqref{ladef2} holds. The complementary slackness condition \eqref{cs1} holds for all $y$ if we take $\gamma_* = b^2$, $\beta_* =0$, and $\alpha_* = (1-b)^2  t$, as before. With $z_* = 0$, the complementary slackness condition \eqref{cs2} holds. The law of $Y$ is primal feasible, and so by Theorem  \ref{thm1} the result follows.

\hfill $\square$

\section{Does more noise mean worse prediction?}\label{S2a}
As we saw at \eqref{normnoise}, if $Y \sim N(0,1)$ is independent of $X$, then
\begin{equation}
E[X|X+\lambda Y] \rightarrow EX \qquad \hbox{\rm a.s.,}
\label{normnoise2}
\end{equation}
so in this situation, {\em adding a larger-variance noise to $X$ decreases the variance of $E[X|X+Y]$.} One might conjecture that this holds more generally, but a little thought shows that this is not so. Indeed, if $X,Y \sim B(1,\half)$, then we have $E[X|X+2Y] = X$, which has larger variance than $E[X|X+Y]$. This being said, a result in the direction of \eqref{normnoise2} is valid if $Y$ is {\em self-decomposable}, as defined in Definition \ref{def1} below.

\begin{defin} \label{def1}
A random variable $Y$ is self-decomposable (belongs to class $\cal L$), if for any $c,\:0<c<1$ there exists a random variable $U_c$ independent of $Y$ such that $Y$ is equal in law to $cY+U_c$.
\end{defin}

All $Y\in{\cal L}$ are infinitely divisible. Not all infinitely divisible random variables are in $\cal L$, but the random variables having stable distributions are in $\cal L$. See Chapter 5 of \cite{Lukacs} for properties of the class $\cal L$.\\
\indent Before stating Theorem \ref{thm2} below, we pause to record a couple of simple facts:
\begin{enumerate}
\item For any random variables $U,\:V,\:W$ with $E|U|<\infty$ and $(U,\:V)$ independent of $W$,
\begin{equation}
E(U|V,\:W)= E(U|V).
\label{fact1}
\end{equation}
\item For any $U,\:V,\:W$ with $E(U^2)<\infty,$
\begin{equation}
{\rm var}E(U|V,\:W)\geq {\rm var}E(U|V+W).
\label{fact2}
\end{equation}
\end{enumerate}
From these, we conclude that if $Z$ is independent of $(X,Y)$ then
\begin{equation}
{\rm var}E(X|X+Y)\geq {\rm var}E(X|X+Y+Z),
\label{needed}
\end{equation}
where we have applied \eqref{fact1} and \eqref{fact2} with $U=X$, $V=X+Y$, $W=Z$.

\begin{theorem}\label{thm2}
Let $X$ be a random variable with ${\rm var}X <\infty$ and $Y\in{\cal L}$. Let $V(\lambda):={\rm var}E(X|X+\lambda Y)$. Then $V(\lambda)$ is monotone decreasing on $(0,\:\infty)$ and monotone increasing on $(-\infty,\:0)$.
\end{theorem}

\medskip\noindent
{\sc Proof}. Let $0<\lambda_1 <\lambda_2$ and set $\lambda_1 =c\lambda_2$ with $0<c<1$. Suppose that $X$, $Y$, $U_c$ are independent random variables with the self-decomposable property
\begin{equation*}
Y \sim c Y + U_c. 
\end{equation*}
Then
\begin{eqnarray*}
\var E[X|X+\lambda_2 Y] &=& \var E[X|X+\lambda_2 cY + \lambda_2 U_c]
\\
&=& \var E[X|X+\lambda_1 Y + \lambda_2 U_c]
\\
&\leq& \var E[X|X+\lambda_1 Y],
\end{eqnarray*}
where the last step follows by \eqref{needed}. Monotonicity in $(-\infty,0)$ follows because $-Y \in \mathcal L$.
\hfill$\square$

%%%%%%%%%%%%%%%%%%%%%%%%%%%%%%%%%%%%%%%%%%%%%%%%%%%%%%%%%%%%%%%%%%%%%%%%%
\section{Examples.}\label{S3}
Our first example is $X \sim B(1,p)$, which is simple enough to allow fairly complete analysis for small $\varepsilon$. Thereafter we take a few examples where $X$ has a symmetric discrete distribution and present numerical solutions. Notice that if $X$ is an integer-valued random variable, and $Y$ is an independent random variable for which $\var E[X|X+Y]$ is minimized, then if we use the integer part $[Y]$ of $Y$ instead, the variance of $E[X|X+Y]$ will be the same. Thus if $X$ is integer-valued, we need only search for minimizing $Y$ among integer-valued $Y$.

\subsection{Binomial distribution.}\label{ss22}
Suppose that $P(X=1) = p = 1-q = P(X=0)$, and that the variance of $Y$ is bounded by $\varepsilon >0 $ as before. We shall assume without loss of generality that $p \geq q$. If $\varepsilon$ is small enough (see \eqref{veps_star}), we conjecture that the optimal $Y$ will take only values 0 and 1, $P[Y=1] = p'=1-q'$, and we will prove that this conjecture is true. As is easily checked,
\begin{equation*}
E[X|X+Y=0] = 0,\quad E[X|X+Y=1] = \frac{pq'}{pq'+qp'},\quad E[X|X+Y=2] = 1.
\end{equation*}
Hence
\begin{equation}
E\bigl( \; E[X|X+Y]^2 \;\bigr) = pp' + \frac{(pq')^2}{pq'+qp'}.  
\label{bin1}
\end{equation}
Some routine calculus shows that this is minimized over $p' \equiv 1-q'$ when
\begin{equation}
p' =  p^* \equiv \frac{\sqrt{p}}{\sqrt{p} + \sqrt{q}},
\label{pdash}
\end{equation}
which is at least $\half$ by our assumption that $p \geq q$.
% The minimized value of \eqref{bin1} is
%\begin{equation}
%\frac{(\sqrt{p}+2\sqrt{q}) p^{3/2}}{(\sqrt{p}+\sqrt{q})^2}
%\label{minvar}
%\end{equation}
%and the variance of $Y$ is 
%\begin{equation}
%\var(Y) 	= \frac{\sqrt{pq}}{(\sqrt{p} + \sqrt{q})^2}.
%\label{varY}
%\end{equation}
Bearing in mind that the variance of $Y$ is bounded by $\varepsilon$, it may be that the value $p^*$ for  $p'$ obtained at \eqref{pdash} is not achievable; if the variance of the $B(1,p')$ random variable $Y$ is bounded by $\varepsilon$, then 
\begin{equation}
 P[Y=1]  \geq \half\{1+\sqrt{1-4\varepsilon}\}.
\label{PYlb}
\end{equation}
We see that this lower bound is greater than or equal to $p^*$ if and only if
\begin{equation}
\varepsilon \leq \quarter (1 - (2p^*-1)^2 ) \equiv \varepsilon_*.
\label{veps_star}
\end{equation}
If $ \varepsilon> \varepsilon_*$, the optimal $Y$ will put positive probability on more than two integer values. {\em  So we proceed on the assumption that $\varepsilon \leq \varepsilon_*$}. With a view to applying Theorem \ref{thm1}, 
we see from \eqref{ladef2} that we must take $\lambda_*$ satisfying
\begin{equation}
\lambda_*(0) =0, \quad \lambda_*(1) = \frac{pq'}{pq' + qp'},
\quad \lambda_*(2) = 1,
\label{lambda_star}
\end{equation}
and now we must check for suitable choices of $\alpha_*$, $\beta_*$, $\gamma_*$, $\lambda$ that the dual-feasibility condition \eqref{df1} holds for  {\em all} integer $y$:
\begin{equation}
q \lambda_*(y)^2 + p (1-\lambda_*(1+y))^2  \leq Q(y) \equiv \alpha_* + \beta_* y + \gamma_* y^2,
\label{DF}
\end{equation}
with equality when $y=0,1$.  We have from \eqref{PYlb}  that $p' = \half (1+\theta)$, where for brevity $\theta \equiv \sqrt{1-4\varepsilon}$. From the fact that \eqref{DF} must hold with equality at $y = 0, 1$, straightforward algebra yields
\begin{equation}
Q(0) = q B, \qquad  Q(1) = p B,
\label{Q01}
\end{equation}
where 
\begin{equation}
B \equiv \frac{pq (1+\theta)^2}{( q(1+\theta) + p(1-\theta))^2}.
\label{Bdef}
\end{equation}
The quadratic 
\begin{equation}
Q_*(x) = B(pq + (x-q)^2) + A x(x-1)
\label{Qstar}
\end{equation}
fits the conditions \eqref{Q01} for any $A \geq 0$. Moreover, if we define 
\begin{equation}
\lambda_*(y) = 0 \quad \forall y <0, \qquad \lambda_*(y) = 1 \quad \forall y > 2,
\end{equation}
we see that by choosing $A>0$ large enough we can ensure \eqref{DF} for all $y \in \Z$ with equality for $y=0,1$. 

\medskip

\begin{figure}[H]
   \hspace*{-2.0cm}\includegraphics[scale=0.32 ]{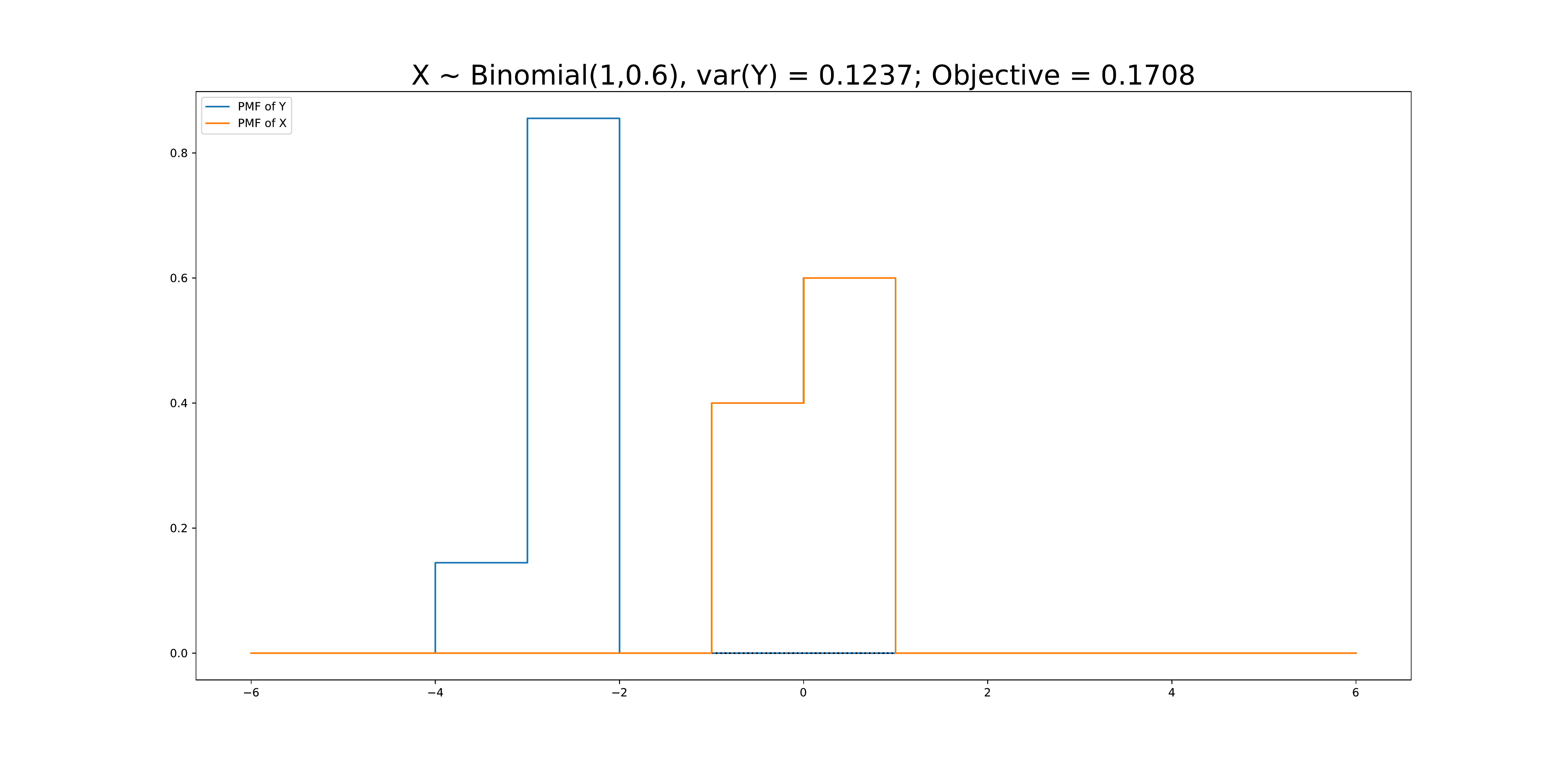} 
    \caption{$X\sim B(1,0.6)$ with low bound on $\var(Y)$.}  
     \label{figbin0}
\end{figure}

\begin{figure}[H]
   \hspace*{-2.0cm}\includegraphics[scale=0.32 ]{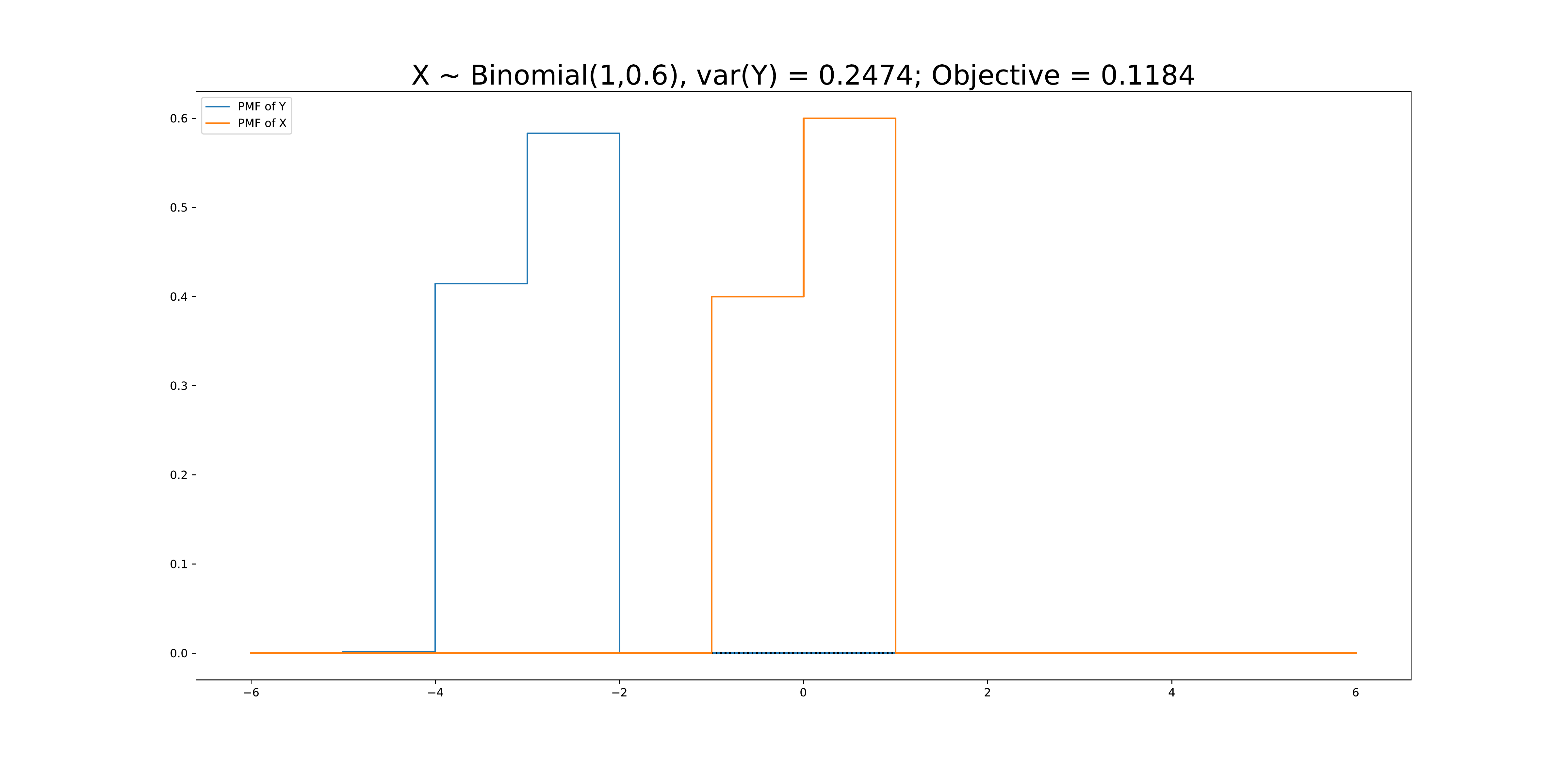} 
    \caption{$X\sim B(1,0.6)$ with critical bound on $\var(Y)$.}  
     \label{figbin1}
\end{figure}

\begin{figure}[H]
   \hspace*{-2.0cm}\includegraphics[scale=0.32 ]{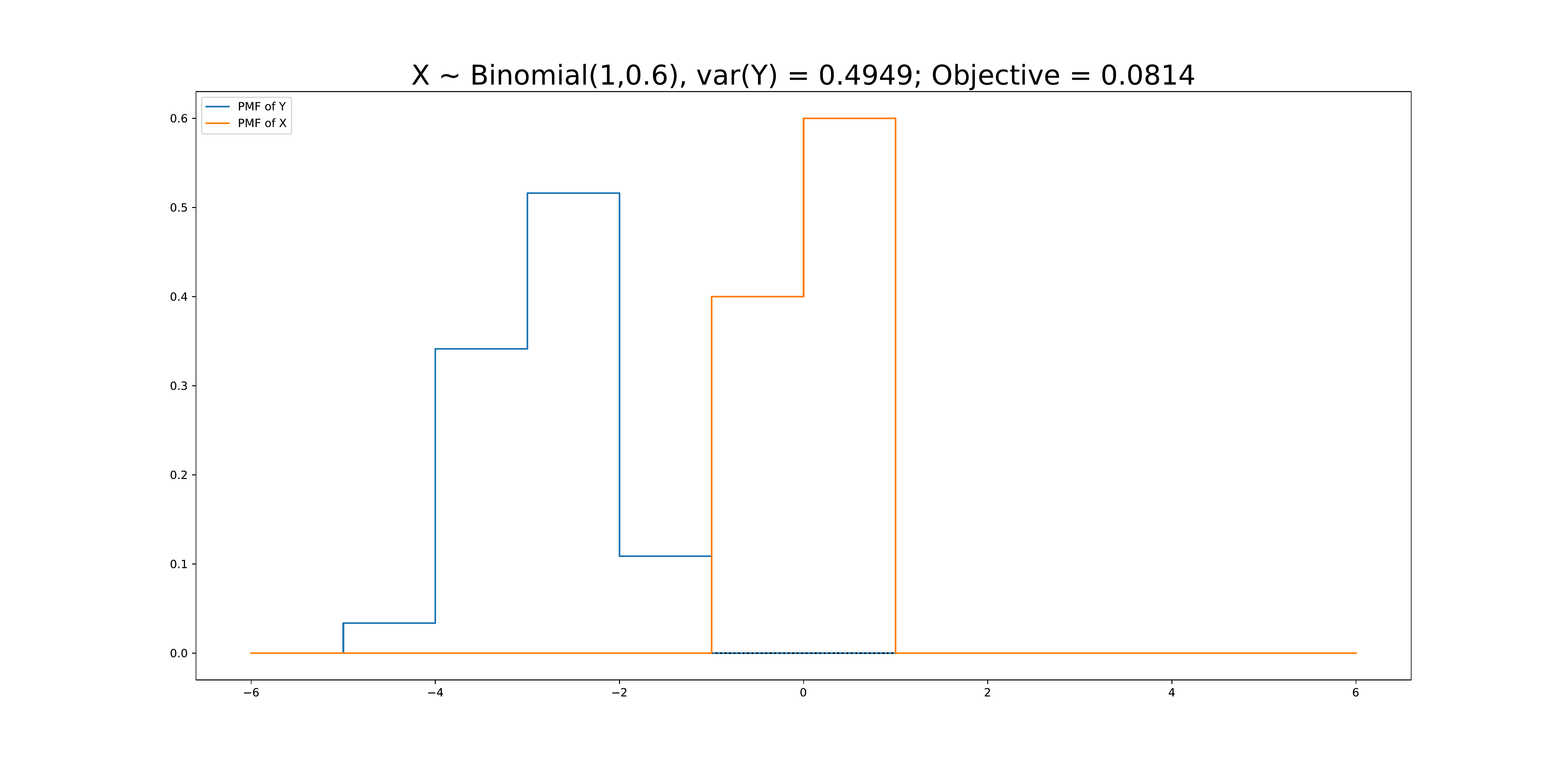} 
    \caption{$X\sim B(1,0.6)$ with high bound on $\var(Y)$.}   
    \label{figbin2}
\end{figure}

To illustrate the kind of solutions we arrive at, we show in Figures \ref{figbin0}, \ref{figbin1} and \ref{figbin2} below the probability mass functions for $X$ and the optimal $Y$ in the case where $X \sim B(1,0.6)$ and $Y$ has to satisfy a low variance bound $\varepsilon = 0.5 \varepsilon_*$, the critical variance bound $\varepsilon_*$, and a higher variance bound $2 \varepsilon_*$ respectively. The probability mass function (PMF) of $Y$ is shown shifted to the left for clarity - as we have already remarked, such a shift makes no difference to the objective. Notice how the objective decreases as the bound on the variance of $Y$ becomes more relaxed, as it should do. Notice also that the PMF of $Y$ in the final plot gives non-zero weight to {\em more than two} values, again as we should expect from the preceding analysis.


%%%%%%%%%%%%%%%%%%%%%%%%%%%%%%%%%%%%%%%%%%%%%%%%%%%%%%%%%%%%%%%%%%%%%%%%%
\subsection{$X$ is uniform.}\label{ss52}
Here we compute the optimal distribution for $Y$ when $X$ is uniform. We consider two cases: the first  low-variance case has $\var(Y) = 2 \var(X)/\pi$ and the second high-variance case has $\var(Y) = \pi \var(X)/2$. The two corresponding figures, Figures \ref{fig2lo} and \ref{fig2hi} below, display the PMFs of $X$ and $Y$, along with a diagnostic plot\footnote{... scaled to fit the plot of the PMFs...} in red and green markers of the computed function
\begin{equation*}
y \mapsto \int f(v) \lambda_*(y-v) \{ \; \lambda_*(y-v) + 2v\;\} \; dv,
\end{equation*}
which according to \eqref{df1} must be dominated by a quadratic\footnote{Recall that $f$ is symmetric.}, and equal to that quadratic wherever the PMF of $Y$ is positive. From our discussion in Section \ref{S2}, if we set $\varepsilon = 2 \var(X)$ then the optimal choice would be to take $Y$ to be the sum of two independent copies of $X$, which in this case would be the sum of two independent uniforms; the resulting PMF would be a symmetric piecewise-linear `tent', and looking at Figure \ref{fig2hi} we something that looks approximately like that.

\begin{figure}[H]
   \hspace*{-2.0cm}\includegraphics[scale=0.3 ]{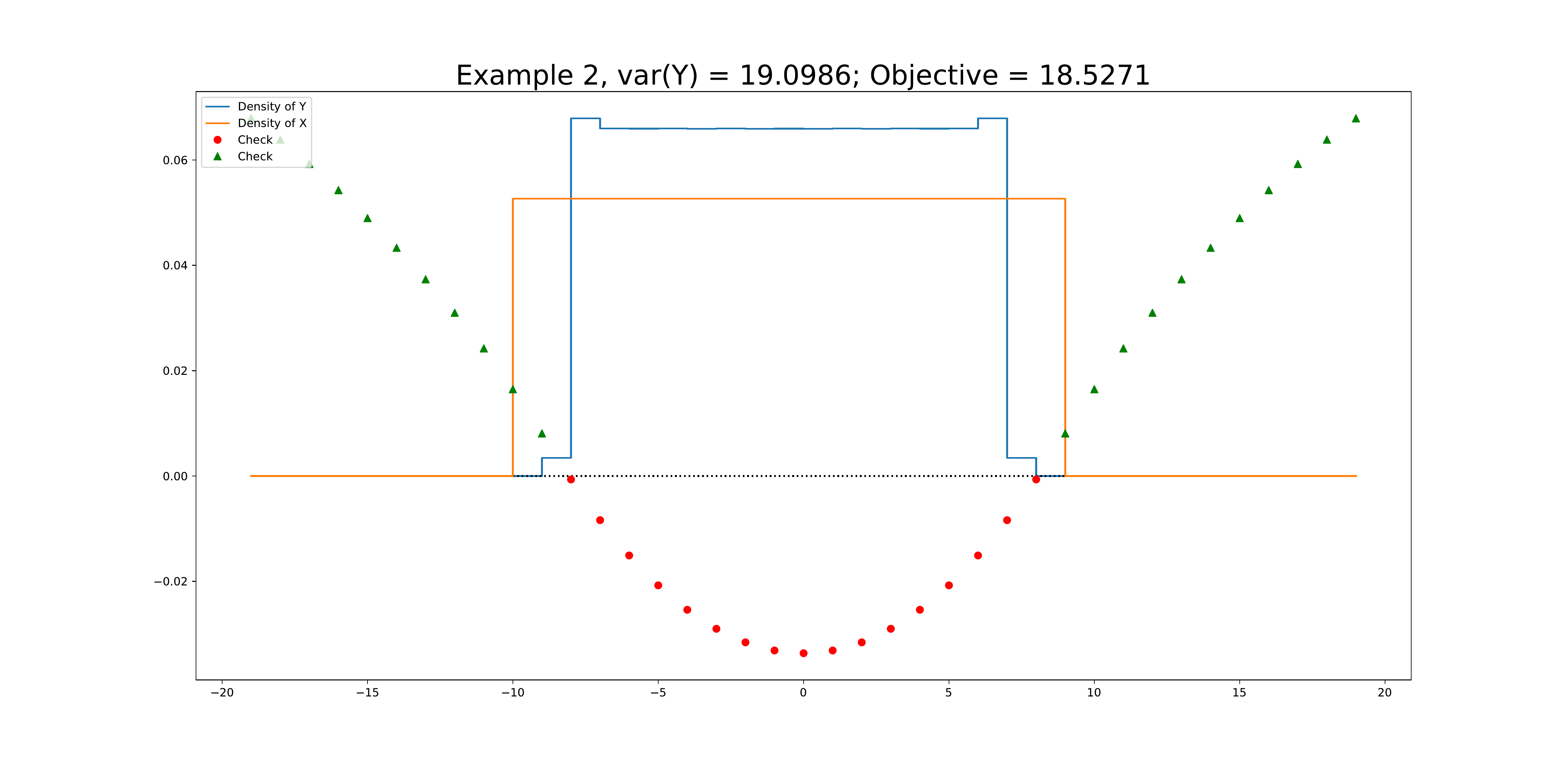} 
    \caption{$X$ uniform, with $\var(Y) = 2 \var(X)/\pi$ .}  
     \label{fig2lo}
\end{figure}
\begin{figure}[H]
   \hspace*{-2.0cm}\includegraphics[scale=0.3 ]{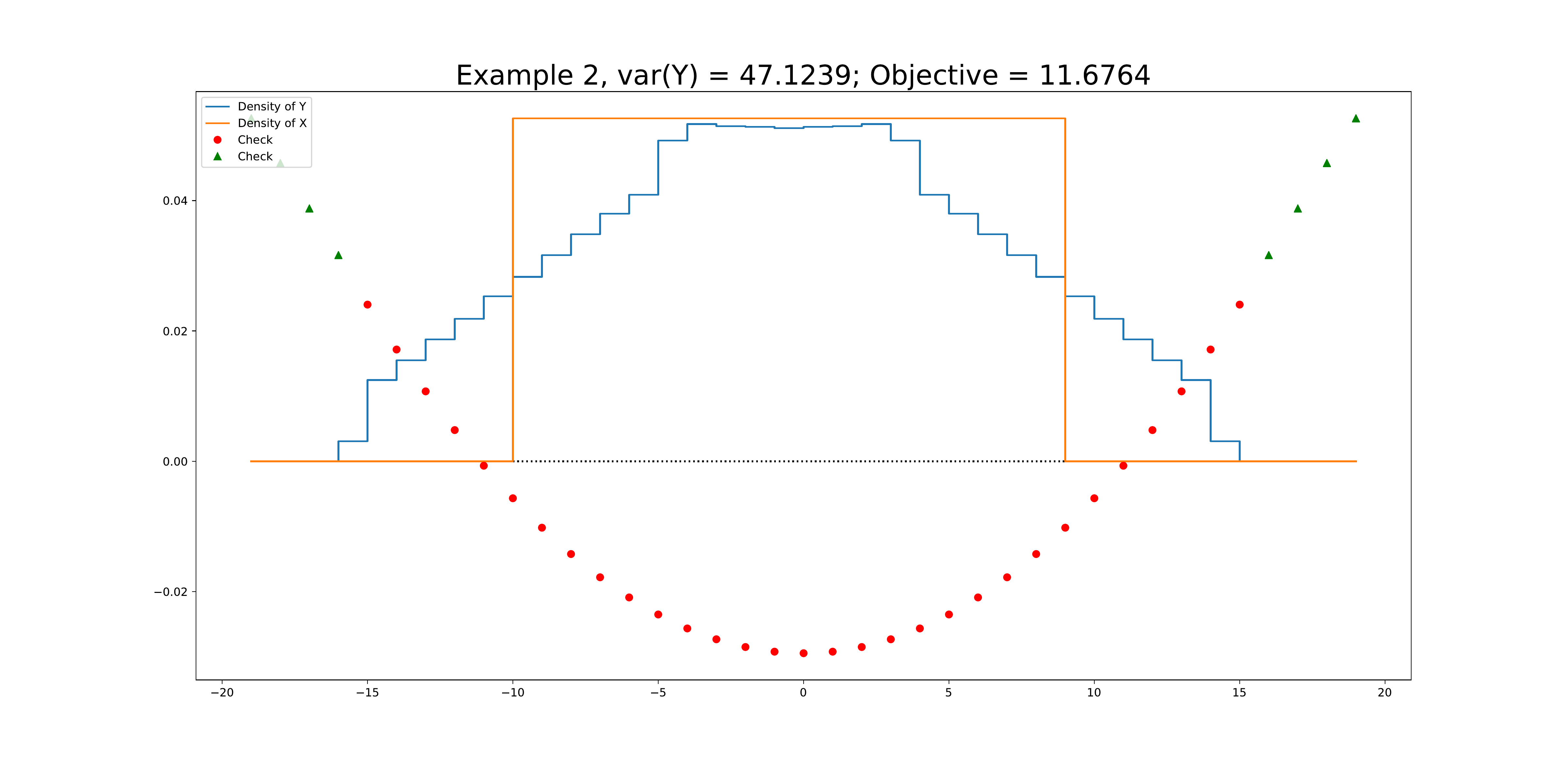} 
    \caption{$X$ uniform, with $\var(Y) = \pi \var(X)/2$.}  
     \label{fig2hi}
\end{figure}

\subsection{$X$ is the sum of two uniforms.}\label{ss53}
Again we compute the optimal $Y$ for two values of $\varepsilon$. Notice how strange the solution is in both cases, particularly for the high variance case, where we see that the distribution of the optimal $Y$ has a hole at the center!
\begin{figure}[H]
   \hspace*{-2.0cm}\includegraphics[scale=0.3 ]{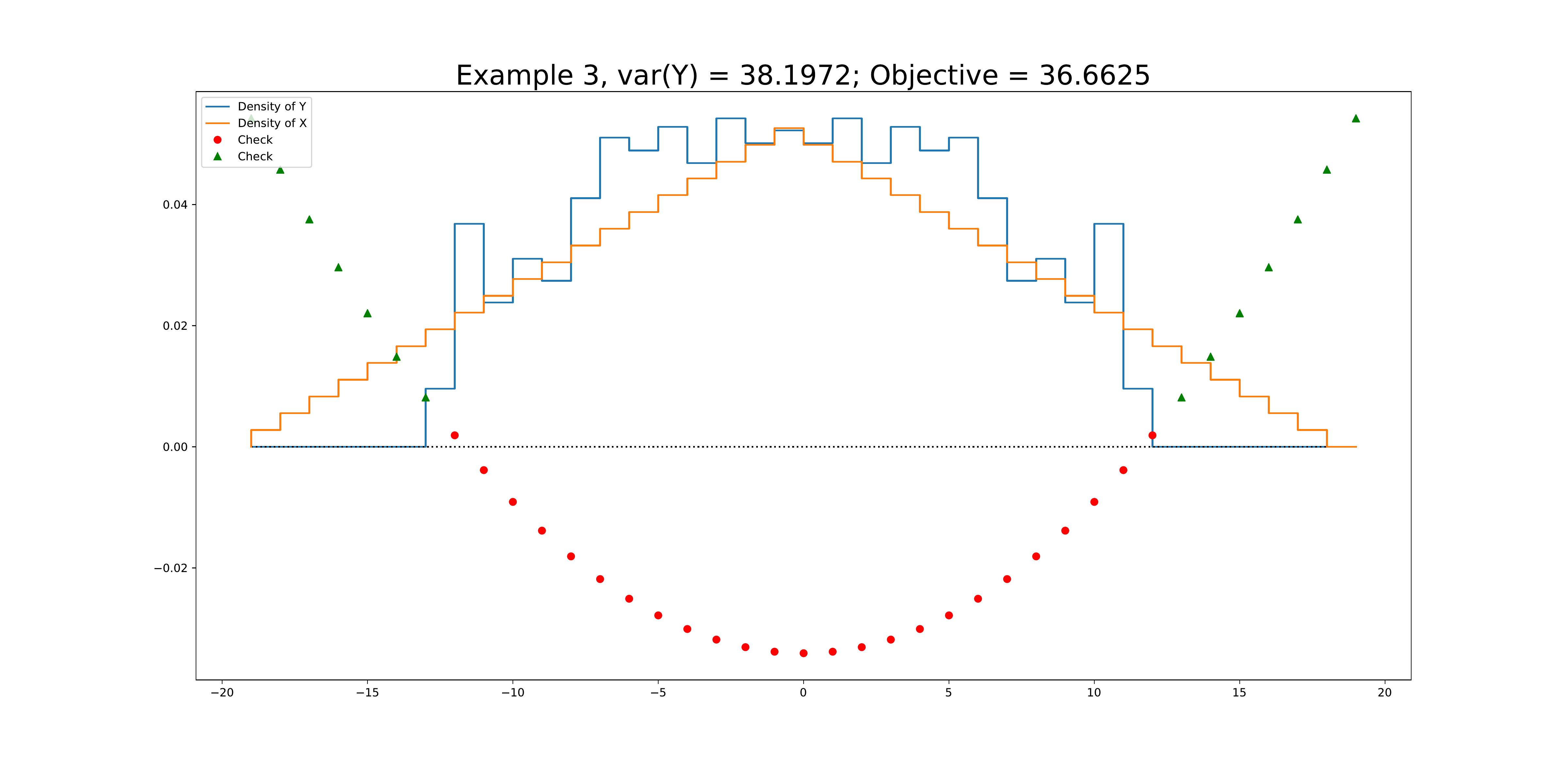} 
    \caption{$X$ is the sum of two uniforms, with $\var(Y) = 2 \var(X)/\pi$ .}  
     \label{fig3lo}
\end{figure}
\begin{figure}[H]
   \hspace*{-2.0cm}\includegraphics[scale=0.3 ]{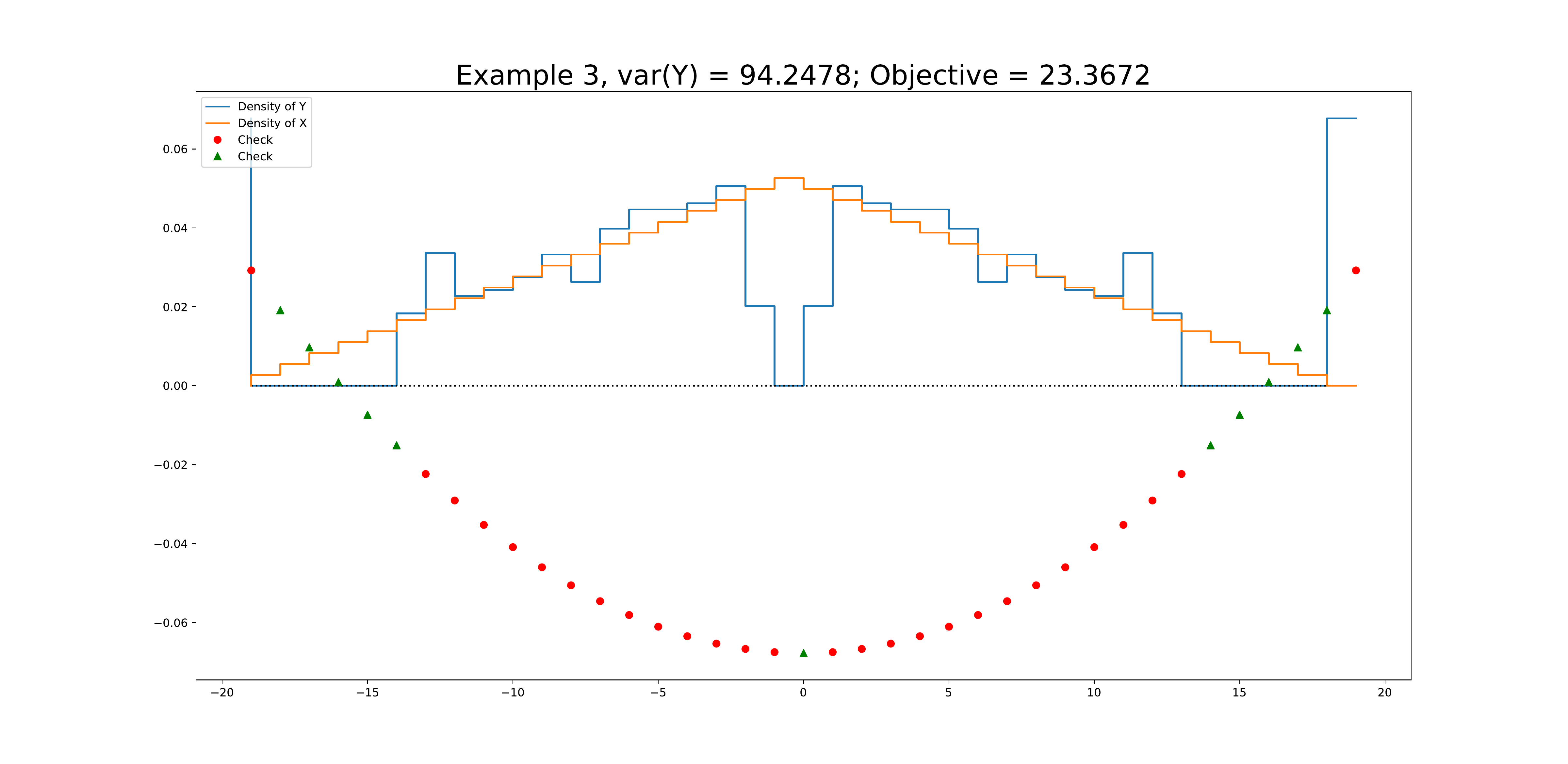} 
    \caption{$X$ is the sum of two uniforms, with $\var(Y) = \pi \var(X)/2$.}  
     \label{fig3hi}
\end{figure}

\subsection{The density of $X$ is the square of that in section \ref{ss53}}\label{ss54}
This time we take the density of $X$ from Section \ref{ss53} and square it (of course, renormalizing to sum to 1).  Once again, the distribution of the optimal $Y$ has a form which would be difficult to guess - the PMF is not monotone in $\Z^+$, for example.

\begin{figure}[H]
   \hspace*{-2.0cm}\includegraphics[scale=0.3 ]{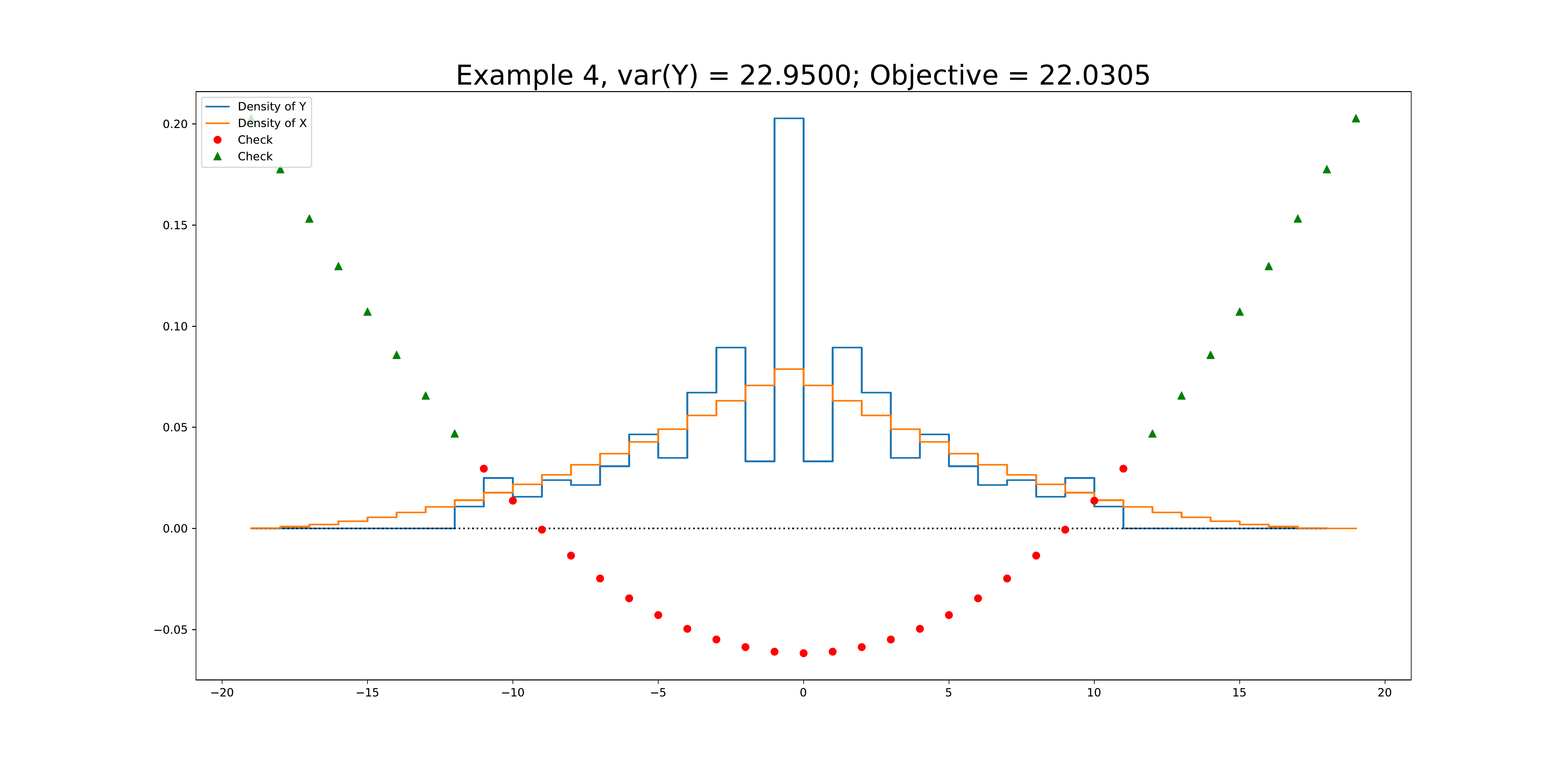} 
    \caption{The density of $X$ is the square of the example in Figure \ref{fig3lo}, with $\var(Y) = 2 \var(X)/\pi$ .}  
     \label{fig4lo}
\end{figure}
\begin{figure}[H]
   \hspace*{-2.0cm}\includegraphics[scale=0.3 ]{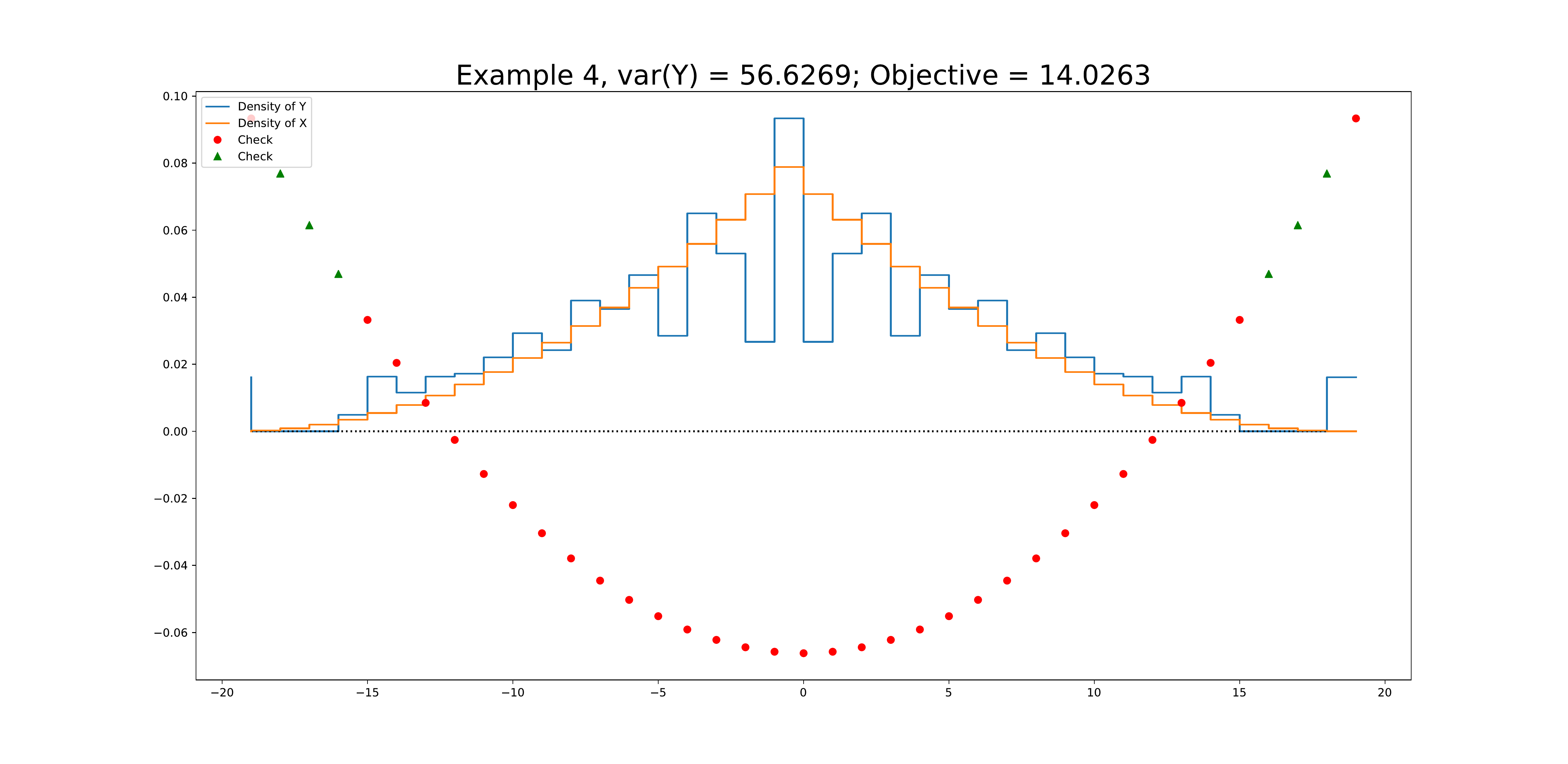} 
    \caption{The density of $X$ is the square of the example in Figure \ref{fig3hi}, with $\var(Y) = \pi \var(X)/2$.}  
     \label{fig4hi}
\end{figure}

\noindent \textbf{Acknowledgments}\\
We thank Professor Dan Crisan and Dongzhou Huang for helpful discussions.

\end{document}